\documentclass[11pt, twoside]{article}
\usepackage{latexsym}
\usepackage{amsmath}
\usepackage{amssymb}
\usepackage[all]{xy}
\usepackage{amsfonts}
\usepackage{verbatim}
\usepackage{amsthm}
\usepackage{mathrsfs}
\usepackage{epsfig}
\usepackage{xy}
\usepackage{array}
\usepackage{stmaryrd}
\usepackage{graphicx,color}
\usepackage{xcolor}
\usepackage{tikz}
\usetikzlibrary{arrows,calc}
\usepackage{etex}
\usepackage{mathdots}
\usepackage{float}
\usepackage{graphics}
\usepackage{pdflscape}
\usepackage{mathrsfs}
\newcommand{\add}{\mathsf{add}\hspace{.01in}}
\usepackage{anysize,hyperref}
\input xypic
\xyoption{all}

\usepackage[perpage,symbol]{footmisc}
\usepackage{setspace}
\renewcommand{\cal}{\mathcal}
\def\A{\mathcal{A}}

\def\B{\mathcal{B}}
\def\C{\mathscr{C}}
\def\E{\mathbb{E}}
\def\s{\mathfrak{s}}

\def\op{^\mathrm{op}}
\def\Ab{\mathit{Ab}}
\def\del{\delta}
\def\dr{\ar@{->}[r]}

\def\X{\mathscr{X}}

\def\add{\mbox{add}}

\def\Hom{\mbox{Hom}}

\begin{document}
\baselineskip=15pt
\title{\Large{\bf Gorenstein dimension of abelian categories}}
\medskip
\author{\textbf{Yu Liu and Panyue Zhou\footnote{Corresponding author. Yu Liu is supported by the Fundamental Research Funds for the Central Universities (Grants No.2682018ZT25). Panyue Zhou was supported by the Hunan Provincial Natural Science Foundation of China (Grant No: 2018JJ3205) and the NSF of China (Grants No.\;11671221)}}}

\date{}

\maketitle
\def\blue{\color{blue}}
\def\red{\color{red}}

\newtheorem{theorem}{Theorem}[section]
\newtheorem{lemma}[theorem]{Lemma}
\newtheorem{corollary}[theorem]{Corollary}
\newtheorem{proposition}[theorem]{Proposition}
\newtheorem{conjecture}{Conjecture}
\theoremstyle{definition}
\newtheorem{definition}[theorem]{Definition}
\newtheorem{question}[theorem]{Question}
\newtheorem{remark}[theorem]{Remark}
\newtheorem{remark*}[]{Remark}
\newtheorem{example}[theorem]{Example}
\newtheorem{example*}[]{Example}
\newtheorem{condition}[theorem]{Condition}
\newtheorem{condition*}[]{Condition}
\newtheorem{construction}[theorem]{Construction}
\newtheorem{construction*}[]{Construction}

\newtheorem{assumption}[theorem]{Assumption}
\newtheorem{assumption*}[]{Assumption}

\baselineskip=17pt
\parindent=0.5cm

\begin{abstract}
\baselineskip=16pt
Let $\C$ be triangulated category and $\X$ a cluster tilting subcategory of $\C$. Koenig and Zhu showed
that the quotient category $\C/\X$ is Gorenstein of Gorenstein dimension at most one.
The notion of an extriangulated category was introduced by Nakaoka and Palu as a simultaneous generalization of exact categories and triangulated categories. Now let $\C$ be extriangulated category with enough projectives and enough injectives, and $\X$ a cluster tilting subcategory of $\C$. In this article, we show that under certain conditions the quotient category $\C/\X$ is Gorenstein of Gorenstein dimension at most one. As an application,
this result generalizes work by Koenig and Zhu. \\[0.5cm]
\textbf{Key words:} Extriangulated categories; abelian quotient categories; Gorenstein dimension.\\[0.2cm]
\textbf{ 2010 Mathematics Subject Classification:} 18E30; 18E10.
\medskip
\end{abstract}

\pagestyle{myheadings}
\markboth{\rightline {\scriptsize   Y. Liu and P. Zhou}}
         {\leftline{\scriptsize  Gorenstein dimension of abelian categories}}

\section{Introduction}

Koenig and Zhu \cite{kz} provided a general framework from  passing from triangulated categories to
abelian categories by factoring out cluster tilting subcategories. More precisely, let $\C$ be a triangulated category and
$\X$ a cluster tilting subcategory of $\C$. They showed that the quotient category $\C/\X$ is an abelian
category and it is Gorenstein of Gorenstein dimension at most one.
Demonet and Liu \cite{dl} gave a way to construct abelian categories from some exact categories.
More precisely, let $\cal B$ be a triangulated category and $\X$ a cluster tilting subcategory of $\cal B$. They showed that the quotient category $\C/\X$ is an abelian
category. Hence, it is quite natural to ask whether this abelian quotient category $\C/\X$ is Gorenstein of Gorenstein dimension at most one. Unfortunately, this result is not always true on an exact category.
See the following example.
\begin{example}\label{ex1}
We revisit Example 3.2 presented in \cite{liu}. Let $\Lambda$ be the $k$-algebra given by the quiver
$$
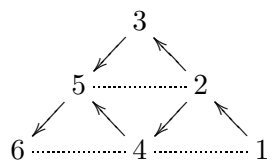
\xymatrix@C=0.4cm@R0.4cm{
&&3 \ar[dl]\\
&5 \ar[dl] \ar@{.}[rr] &&2 \ar[dl] \ar[ul]\\
6 \ar@{.}[rr] &&4 \ar[ul] \ar@{.}[rr] &&1 \ar[ul]}$$
with mesh relations, where $k$ is a field. The AR-quiver of $\B:={\rm mod}\Lambda$ is given by
$$
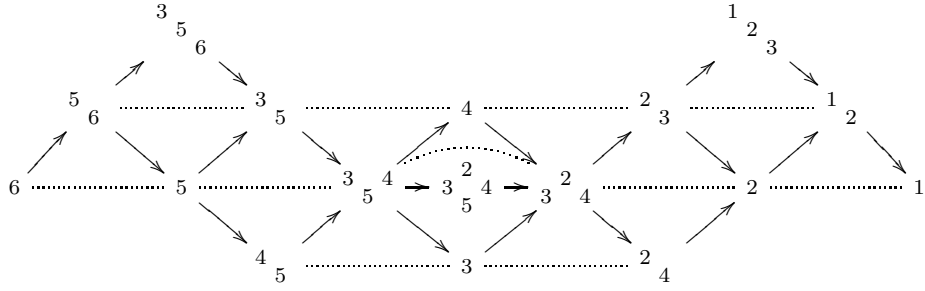
\xymatrix@C=0.3cm@R0.3cm{
&&{\begin{smallmatrix}
3&&\\
&5&\\
&&6
\end{smallmatrix}} \ar[dr] &&&&&&{\begin{smallmatrix}
1&&\\
&2&\\
&&3
\end{smallmatrix}} \ar[dr]\\
&{\begin{smallmatrix}
5&&\\
&6&
\end{smallmatrix}} \ar[ur] \ar@{.}[rr] \ar[dr] &&{\begin{smallmatrix}
3&&\\
&5&
\end{smallmatrix}} \ar@{.}[rr] \ar[dr] &&{\begin{smallmatrix}
4
\end{smallmatrix}} \ar@{.}[rr] \ar[dr] &&{\begin{smallmatrix}
2&&\\
&3&
\end{smallmatrix}} \ar[ur] \ar@{.}[rr] \ar[dr] &&{\begin{smallmatrix}
1&&\\
&2&
\end{smallmatrix}} \ar[dr]\\
{\begin{smallmatrix}
6
\end{smallmatrix}} \ar[ur] \ar@{.}[rr] &&{\begin{smallmatrix}
5
\end{smallmatrix}} \ar[ur] \ar@{.}[rr] \ar[dr] &&{\begin{smallmatrix}
3&&4\\
&5&
\end{smallmatrix}} \ar[ur] \ar[r] \ar[dr] \ar@{.}@/^15pt/[rr] &{\begin{smallmatrix}
&2&\\
3&&4\\
&5&
\end{smallmatrix}} \ar[r] &{\begin{smallmatrix}
&2&\\
3&&4
\end{smallmatrix}} \ar[ur] \ar@{.}[rr] \ar[dr] &&{\begin{smallmatrix}
2
\end{smallmatrix}} \ar[ur] \ar@{.}[rr] &&{\begin{smallmatrix}
1
\end{smallmatrix}}\\
&&&{\begin{smallmatrix}
4&&\\
&5&
\end{smallmatrix}} \ar[ur] \ar@{.}[rr] &&{\begin{smallmatrix}
3
\end{smallmatrix}} \ar[ur] \ar@{.}[rr] &&{\begin{smallmatrix}
2&&\\
&4&
\end{smallmatrix}} \ar[ur]
}$$
We denote by ``$\circ$" in the AR-quiver the indecomposable objects belong to a subcategory. Put
$$
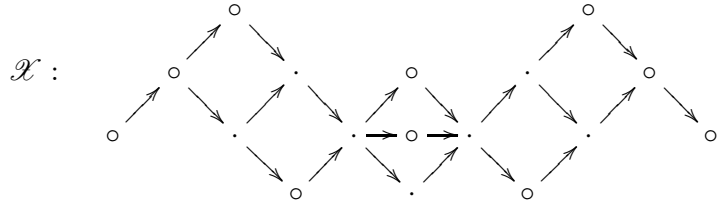
\xymatrix@C=0.4cm@R0.4cm{
&&&\circ \ar[dr] &&&&&&\circ \ar[dr]\\
{\X:} &&\circ \ar[ur]  \ar[dr] &&\cdot  \ar[dr] &&\circ  \ar[dr] &&\cdot  \ar[ur]  \ar[dr] &&\circ \ar[dr]\\
&\circ \ar[ur]  &&\cdot \ar[ur]  \ar[dr] &&\cdot \ar[ur] \ar[r] \ar[dr] &\circ \ar[r] &\cdot \ar[ur] \ar[dr] &&\cdot \ar[ur] &&\circ\\
&&&&\circ \ar[ur] &&\cdot \ar[ur] &&\circ \ar[ur]
}
$$
$\X$ is a cluster tilting subcategory of $\B$. Then $\B/\X\simeq {\rm mod}(\Omega \X/\cal P)$ where $\cal P$ is the subcategory of projective objects and $\Omega$ is a \emph{syzygy} functor, and its quiver is the following:
 $$
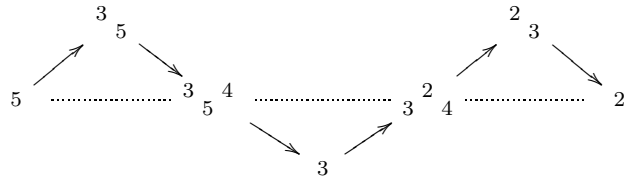
\xymatrix@C=0.4cm@R0.4cm{
&{\begin{smallmatrix}
3&\ \\
&5
\end{smallmatrix}} \ar[dr]
&&&&{\begin{smallmatrix}
2&\ \\
&3
\end{smallmatrix}}\ar[dr]\\
{\begin{smallmatrix}
\ &5&\
\end{smallmatrix}} \ar@{.}[rr] \ar[ur]
&&{\begin{smallmatrix}
3&&4\ \\
&5&
\end{smallmatrix}}\ar[dr] \ar@{.}[rr]
&&{\begin{smallmatrix}
&2&\ \\
3&&4
\end{smallmatrix}}\ar[ur] \ar@{.}[rr]
&&{\begin{smallmatrix}
\ &2&\
\end{smallmatrix}}\\
&&&{\begin{smallmatrix}
\ &3&\
\end{smallmatrix}} \ar[ur]}$$
It is not Gorenstein of Gorenstein dimension at most one. Note that the non-projective injective object $2$ have projective dimension $3$.
\end{example}

Recently, the notion of an extriangulated category was introduced by Nakaoka and Palu \cite{np} as a simultaneous generalization of exact categories and triangulated categories.
Cluster tilting theory gives a way to construct abelian categories from some extriangulated categories.
Let $\C$ be extriangulated category with enough projectives and enough injectives, and $\X$ a cluster tilting subcategory of $\C$. Then the quotient category $\C/\X$ is an abelian category, see \cite{zz,ln}.
We know that a module category can be viewed as an extriangulated category with enough projectives and enough injectives. Hence the abelian quotient category $\C/\X$ is not Gorenstein of Gorenstein dimension at most one,
see Example \ref{ex1}.

Let $\C$ be extriangulated category with enough projectives $\cal P$ and enough injectives $\cal I$, and $\X$ a subcategory of $\C$.
We denote $\Omega\X={\rm CoCone}(\cal P,\X)$, that is to say, $\Omega\X$ is the subcategory of $\C$ consisting of
objects $\Omega X$ such that there exists an $\E$-triangle:
$$\Omega X\overset{a}{\longrightarrow}P\overset{b}{\longrightarrow}X\overset{}{\dashrightarrow},$$
with $P\in\cal P$ and $X\in\X$. We call $\Omega$ the \emph{syzygy} of $\X$.
Dually we define the \emph{cosyzygy} of $\X$ by $\Sigma\X={\rm Cone}(\X,\cal I)$.
Namely, $\Sigma\X$ is the subcategory of $\C$ consisting of
objects $\Sigma X$ such that there exists an $\E$-triangle:
$$X\overset{c}{\longrightarrow}I\overset{d}{\longrightarrow}\Sigma X\overset{}{\dashrightarrow},$$
with $I\in\cal I$ and $X\in\X$. For more details, see \cite[Definition 4.2]{ln} and \cite[Proposition 4.3]{ln}.

Our main result is as follows, which gives sufficient conditions on the quotient category $\C/\X$ is Gorenstein of Gorenstein dimension at most one, where $\C$ is an extriangulated category with enough projectives and enough injectives and $\X$ is a cluster tilting subcategory of $\C$.
\begin{theorem}\emph{(see Theorem \ref{thm} for more details)}\label{thm0}
Let $\C$ be an extriangulated category with enough projective objects  and enough injective objects.
Suppose that $\X$ is a cluster-tilting subcategory of $\C$ and $\A$ is the abelian quotient category $\C/\X$. Then:
\begin{itemize}
\item[\emph{(1)}] The category $\A$ has enough projective objects and enough injective objects.
\item[\emph{(2)}] If $\Sigma(\Omega\X)\subseteq\X$ and $\Omega(\Sigma\X)\subseteq\X$, then the category $\A$ is Gorenstein of Gorenstein dimension at most one.
\end{itemize}
\end{theorem}

We know that a triangulated category can be viewed as extriangulated category with enough projective objects and enough injective objects. In Theorem \ref{thm0}, this condition $\Sigma(\Omega\X)\subseteq\X$ and $\Omega(\Sigma\X)\subseteq\X$ is automatically satisfied. As an application, our result generalizes work by Koenig and Zhu \cite[Theorem 4.3]{kz}.

The paper is organised as follows: In Section 2, we review some elementary definitions and facts
that we need to use. In Section 3, we prove the main result of this article.

\section{Preliminaries}

We recall some definitions and basic properties of extriangulated categories from \cite{np}.
Let $\C$ be an additive category. Suppose that $\C$ is equipped with a biadditive functor $$\E\colon\C\op\times\C\to\Ab,$$
where $\Ab$ is the category of abelian groups. For any pair of objects $A,C\in\C$, an element $\delta\in\E(C,A)$ is called an {\it $\E$-extension}. Thus formally, an $\E$-extension is a triplet $(A,\delta,C)$.
Let $(A,\del,C)$ be an $\E$-extension. Since $\E$ is a bifunctor, for any $a\in\C(A,A')$ and $c\in\C(C',C)$, we have $\E$-extensions
$$ \E(C,a)(\del)\in\E(C,A')\ \ \text{and}\ \ \ \E(c,A)(\del)\in\E(C',A). $$
We abbreviately denote them by $a_\ast\del$ and $c^\ast\del$.
For any $A,C\in\C$, the zero element $0\in\E(C,A)$ is called the \emph{spilt $\E$-extension}.

\begin{definition}{\cite[Definition 2.3]{np}}
Let $(A,\del,C),(A',\del',C')$ be any pair of $\E$-extensions. A {\it morphism} $$(a,c)\colon(A,\del,C)\to(A',\del',C')$$ of $\E$-extensions is a pair of morphisms $a\in\C(A,A')$ and $c\in\C(C,C')$ in $\C$, satisfying the equality
$ a_\ast\del=c^\ast\del'. $
Simply we denote it as $(a,c)\colon\del\to\del'$.
\end{definition}

\begin{definition}{\cite[Definition 2.6]{np}}
Let $\delta=(A,\delta,C),\delta^{\prime}=(A^{\prime},\delta^{\prime},C^{\prime})$ be any pair of $\mathbb{E}$-extensions. Let
\[ C\overset{\iota_C}{\longrightarrow}C\oplus C^{\prime}\overset{\iota_{C^{\prime}}}{\longleftarrow}C^{\prime} \]
and
\[ A\overset{p_A}{\longleftarrow}A\oplus A^{\prime}\overset{p_{A^{\prime}}}{\longrightarrow}A^{\prime} \]
be coproduct and product in $\B$, respectively. Since $\mathbb{E}$ is biadditive, we have a natural isomorphism
\[ \mathbb{E}(C\oplus C^{\prime},A\oplus A^{\prime})\cong \mathbb{E}(C,A)\oplus\mathbb{E}(C,A^{\prime})\oplus\mathbb{E}(C^{\prime},A)\oplus\mathbb{E}(C^{\prime},A^{\prime}). \]

Let $\delta\oplus\delta^{\prime}\in\mathbb{E}(C\oplus C^{\prime},A\oplus A^{\prime})$ be the element corresponding to $(\delta,0,0,\delta^{\prime})$ through the above isomorphism. This is the unique element which satisfies
\begin{eqnarray*}
\mathbb{E}(\iota_C,p_A)(\delta\oplus\delta^{\prime})=\delta&,&\mathbb{E}(\iota_C,p_{A^{\prime}})(\delta\oplus\delta^{\prime})=0,\\
\mathbb{E}(\iota_{C^{\prime}},p_A)(\delta\oplus\delta^{\prime})=0&,&\mathbb{E}(\iota_{C^{\prime}},p_{A^{\prime}})(\delta\oplus\delta^{\prime})=\delta^{\prime}.
\end{eqnarray*}
\end{definition}

\begin{definition}{\cite[Definition 2.7]{np} and \cite[Definition 2.8]{np}}
Let $A,C\in\C$ be any pair of objects. Sequences of morphisms in $\C$
$$\xymatrix@C=0.7cm{A\ar[r]^{x} & B \ar[r]^{y} & C}\ \ \text{and}\ \ \ \xymatrix@C=0.7cm{A\ar[r]^{x'} & B' \ar[r]^{y'} & C}$$
are said to be {\it equivalent} if there exists an isomorphism $b\in\C(B,B')$ which makes the following diagram commutative.
$$\xymatrix{
A \ar[r]^x \ar@{=}[d] & B\ar[r]^y \ar[d]_{\simeq}^{b} & C\ar@{=}[d]&\\
A\ar[r]^{x'} & B' \ar[r]^{y'} & C &}$$

We denote the equivalence class of $\xymatrix@C=0.7cm{A\ar[r]^{x} & B \ar[r]^{y} & C}$ by $[\xymatrix@C=0.7cm{A\ar[r]^{x} & B \ar[r]^{y} & C}]$.

For any $A,C\in\C$, we denote $0=[
A \xrightarrow{\binom{1}{0}} A\oplus C \xrightarrow{(0,1)} C ]$.

For any two equivalence classes, we denote
$ [A\overset{x}{\longrightarrow}B\overset{y}{\longrightarrow}C]\oplus [A^{\prime}\overset{x^{\prime}}{\longrightarrow}B^{\prime}\overset{y^{\prime}}{\longrightarrow}C^{\prime}]=[A\oplus A^{\prime}\overset{x\oplus x^{\prime}}{\longrightarrow}B\oplus B^{\prime}\overset{y\oplus y^{\prime}}{\longrightarrow}C\oplus C^{\prime}]. $
\end{definition}

\begin{definition}{\cite[Definition 2.9]{np}}
Let $\s$ be a correspondence which associates an equivalence class $\s(\del)=[\xymatrix@C=0.7cm{A\ar[r]^{x} & B \ar[r]^{y} & C}]$ to any $\E$-extension $\del\in\E(C,A)$. This $\s$ is called a {\it realization} of $\E$, if it satisfies the following condition:
\begin{itemize}
\item Let $\del\in\E(C,A)$ and $\del'\in\E(C',A')$ be any pair of $\E$-extensions, with $$\s(\del)=[\xymatrix@C=0.7cm{A\ar[r]^{x} & B \ar[r]^{y} & C}],\ \ \ \s(\del')=[\xymatrix@C=0.7cm{A'\ar[r]^{x'} & B'\ar[r]^{y'} & C'}].$$
Then, for any morphism $(a,c)\colon\del\to\del'$, there exists $b\in\C(B,B')$ which makes the following diagram commutative.
$$\xymatrix{
A \ar[r]^x \ar[d]^a & B\ar[r]^y \ar[d]^{b} & C\ar[d]^c&\\
A'\ar[r]^{x'} & B' \ar[r]^{y'} & C' &}$$
\end{itemize}
In the above situation, we say that the triplet $(a,b,c)$ realizes $(a,c)$.
\end{definition}

\begin{definition}{\cite[Definition 2.10]{np}}
A realization $\s$ of $\E$ is called \emph{additive} if it satisfies the following conditions.
\begin{itemize}
\item[(1)] For any $A,C\in\C$, the split $\E$-extension $0\in\E(C,A)$ satisfies $\s(0)=0$.
\item[(2)] For any pair of $\E$-extensions $\delta\in\E(C,A)$ and $\delta'\in\E(C',A')$,
$$\s(\delta\oplus\delta')=\s(\delta)\oplus\s(\delta')$$
holds.
\end{itemize}
\end{definition}

\begin{definition}{\cite[Definition 2.12]{np}}
A triplet $(\C,\E,\s)$ is called an \emph{externally triangulated category} (or \emph{extriangulated category} for short) if it satisfies the following conditions:
\begin{itemize}
\item[{\rm (ET1)}] $\E\colon\C\op\times\C\to\Ab$ is a biadditive functor.
\item[{\rm (ET2)}] $\s$ is an additive realization of $\E$.
\item[{\rm (ET3)}] Let $\del\in\E(C,A)$ and $\del'\in\E(C',A')$ be any pair of $\E$-extensions, realized as
$$ \s(\del)=[\xymatrix@C=0.7cm{A\ar[r]^{x} & B \ar[r]^{y} & C}],\ \ \s(\del')=[\xymatrix@C=0.7cm{A'\ar[r]^{x'} & B' \ar[r]^{y'} & C'}]. $$
For any commutative square
$$\xymatrix{
A \ar[r]^x \ar[d]^a & B\ar[r]^y \ar[d]^{b} & C&\\
A'\ar[r]^{x'} & B' \ar[r]^{y'} & C' &}$$
in $\C$, there exists a morphism $(a,c)\colon\del\to\del'$ satisfying $cy=y'b$.
\item[{\rm (ET3)$\op$}] Dual of (ET3).
\item[{\rm (ET4)}] Let $(A,\del,D)$ and $(B,\del',F)$ be $\E$-extensions realized by
$$\xymatrix@C=0.7cm{A\ar[r]^{f} & B \ar[r]^{f'} & D}\ \ \text{and}\ \ \ \xymatrix@C=0.7cm{B\ar[r]^{g} & C \ar[r]^{g'} & F}$$
respectively. Then there exist an object $E\in\C$, a commutative diagram
$$\xymatrix{A\ar[r]^{f}\ar@{=}[d]&B\ar[r]^{f'}\ar[d]^{g}&D\ar[d]^{d}\\
A\ar[r]^{h}&C\ar[d]^{g'}\ar[r]^{h'}&E\ar[d]^{e}\\
&F\ar@{=}[r]&F}$$
in $\C$, and an $\E$-extension $\del^{''}\in\E(E,A)$ realized by $\xymatrix@C=0.7cm{A\ar[r]^{h} & C \ar[r]^{h'} & E},$ which satisfy the following compatibilities.
\begin{itemize}
\item[{\rm (i)}] $\xymatrix@C=0.7cm{D\ar[r]^{d} & E \ar[r]^{e} & F}$  realizes $f'_{\ast}\del'$,
\item[{\rm (ii)}] $d^\ast\del''=\del$,

\item[{\rm (iii)}] $f_{\ast}\del''=e^{\ast}\del'$.
\end{itemize}

\item[{\rm (ET4)$\op$}]  Dual of (ET4).
\end{itemize}
\end{definition}

\begin{remark}
We know that both exact categories and triangulated categories are extriangulated categories
\cite[Example 2.13]{np} and extension-closed subcategories of extriangulated categories are
again extriangulated \cite[Remark 2.18]{np}. Moreover, there exist extriangulated categories which
are neither exact categories nor triangulated categories \cite[Proposition 3.30]{np} and \cite[Example 4.14]{zz1}.
\end{remark}

We use the following terminology.
\begin{definition}{\cite{np}}
Let $(\C,\E,\s)$ be an extriangulated category.
\begin{itemize}
\item[(1)] If a conflation $A\overset{x}{\longrightarrow}B\overset{y}{\longrightarrow}C$ realizes $\delta\in\mathbb{E}(C,A)$, we call the pair $(A\overset{x}{\longrightarrow}B\overset{y}{\longrightarrow}C,\delta)$ an {\it $\E$-triangle}, and write it in the following way.
$$A\overset{x}{\longrightarrow}B\overset{y}{\longrightarrow}C\overset{\delta}{\dashrightarrow}$$

\item[(2)] Let $A\overset{x}{\longrightarrow}B\overset{y}{\longrightarrow}C\overset{\delta}{\dashrightarrow}$ and $A^{\prime}\overset{x^{\prime}}{\longrightarrow}B^{\prime}\overset{y^{\prime}}{\longrightarrow}C^{\prime}\overset{\delta^{\prime}}{\dashrightarrow}$ be any pair of $\E$-triangles. If a triplet $(a,b,c)$ realizes $(a,c)\colon\delta\to\delta^{\prime}$, then we write it as
$$\xymatrix{
A \ar[r]^x \ar[d]^a & B\ar[r]^y \ar[d]^{b} & C\ar@{-->}[r]^{\del}\ar[d]^c&\\
A'\ar[r]^{x'} & B' \ar[r]^{y'} & C'\ar@{-->}[r]^{\del'} &}$$
and call $(a,b,c)$ a {\it morphism of $\E$-triangles}.

\item[(3)] An object $P\in\C$ is called {\it projective} if
for any $\E$-triangle $A\overset{x}{\longrightarrow}B\overset{y}{\longrightarrow}C\overset{\delta}{\dashrightarrow}$ and any morphism $c\in\C(P,C)$, there exists $b\in\C(P,B)$ satisfying $yb=c$.
We denote the subcategory of projective objects by $\cal P\subseteq\C$. Dually, the subcategory of injective objects is denoted by $\cal I\subseteq\C$.

\item[(4)] We say that $\C$ {\it has enough projective objects} if
for any object $C\in\C$, there exists an $\E$-triangle
$A\overset{x}{\longrightarrow}P\overset{y}{\longrightarrow}C\overset{\delta}{\dashrightarrow}$
satisfying $P\in\cal P$. We can define the notion of having enough injectives dually.

\end{itemize}
\end{definition}

\begin{definition}\cite[Definition 2.10]{zz}
Let $(\C, \E, \s)$ be an extriangulated category and $\X$ a subcategory of $\C$.
\begin{itemize}
\item $\X$ is called \emph{rigid} if $\E(\X,\X)=0$;
\item $\X$ is called \emph{cluster tilting} if it satisfies the following conditions:
\begin{enumerate}
\item[(1)] $\X$ is a functorially finite in $\C$;
\item[(2)] $M\in \X$ if and only if $\E(M,\X)=0$;
\item[(3)] $M\in \X$ if and only if $\E(\X,M)=0$.
\end{enumerate}
\end{itemize}
\end{definition}

By definition of a cluster tilting subcategory, we can immediately conclude:
\begin{remark}\label{y6}
Let $(\C, \E, \s)$ be an extriangulated category with enough projectives and enough injectives.
\begin{itemize}
\item  If $\X$ is a cluster tilting  subcategory of $\C$, then $\cal P\subseteq\X$ and $\cal I\subseteq\X$.
\item  $\X$ is a cluster tilting subcategory of $\C$ if and only if
\begin{enumerate}
\item[(1)] $\X$ is rigid;
\item[(2)] For any $C\in\C$, there exists an $\E$-triangle $\xymatrix@C=0.5cm{C\ar[r]^{a\;} & X_1 \ar[r]^{b} & X_2\ar@{-->}[r]^{\del}&}$, where $X_1,X_2\in\X$;
\item[(3)] For any $C\in\C$, there exists an $\E$-triangle $\xymatrix@C=0.5cm{X_3\ar[r]^{c\;} & X_4 \ar[r]^{d} & C\ar@{-->}[r]^{\eta}&}$, where $X_3,X_4\in\X$.
\end{enumerate}
\end{itemize}
\end{remark}

\begin{theorem}\emph{\cite[Theorem 3.4]{zz} and \cite[Theorem 3.2]{ln}}
let $\C$ be extriangulated category with enough projectives and enough injectives, and $\X$ a cluster tilting subcategory of $\C$. The quotient category $\C/\X$ is an abelian category.
\end{theorem}

\section{Gorenstein dimension at most one}
Let $\C$ be an additive category and $\X$ be a subcategory of $\C$. We denote by $\C/\X$
the category whose objects are objects of $\C$ and whose morphisms are elements of
$\Hom_{\C}(A,B)/\X(A,B)$ for $A,B\in\C$, where $\X(A,B)$ the subgroup of $\Hom_{\C}(A,B)$ consisting of morphisms
which factor through an object in $\X$.
Such category is called the quotient category
of $\C$ by $\X$. For any morphism $f\colon A\to B$ in $\C$, we denote by $\overline{f}$ the image of $f$ under
the natural quotient functor $\C\to\C/\X$.

\begin{lemma}\label{x2}
Let $\C$ be an extriangulated category with enough projectives  and enough injectives.
Suppose that $\X$ is a cluster-tilting subcategory of $\C$ and $\A$ is the abelian quotient category $\C/\X$.
Then any object $C\in\C$ admits an epimorphism $\underline{\beta}\colon \Omega X\to C$ for some $X\in\X$ in $\A$. Dually any object $C\in\C$ admits a monomorphism $\underline{\alpha}\colon C\to \Sigma X$ for some $X\in\X$ in $\A$.
\end{lemma}

\proof We only prove the first statement. The second statement is dual.

Since $\X$ is cluster tilting, there exists an $\E$-triangle
$$C\xrightarrow{~a~}X_0\xrightarrow{~b~}X_1\overset{}{\dashrightarrow},$$
where $X_0,X_1\in\X$. By definition $\Omega X_0$ admits an $\E$-triangle
$$\Omega X_0\xrightarrow{~u~}P\xrightarrow{~v~}X_0\overset{}{\dashrightarrow},$$
By (ET4)$\op$, we have the following commutative diagram made of $\E$-triangles
\begin{equation}\label{t0}
\begin{array}{l}
\xymatrix{\Omega X_0\ar[d]\ar@{=}[r]&\Omega X_0\ar[d]&&\\
\Omega X_1\ar[r]^u\ar[d]^{\beta}&P\ar[r]^v\ar[d]^{\gamma}&X_1\ar@{=}[d]\\
C\ar[r]^{a}&X_0\ar[r]^{b}&X_1}
\end{array}
\end{equation}
We claim that $\overline{\beta}\colon \Omega X_1\to C$ is an epimorphism in $\A$.
In fact, assume that $\overline{c}\colon C\to B$ is any morphism in $\A$ such that
$\overline{c}\circ\overline{\beta}=0$. Then $c\beta$ factors through $\X$. Since $u$
is a left $\X$-approximation of $\Omega X_1$, there exists a morphism
$w\colon P\to B$ such that $c\beta=wu$.

By \cite[Lemma 3.13]{np}, the lower-left square in the diagram (\ref{t0})
$$\xymatrix{\Omega X_1\ar[r]^{u}\ar[d]^{\beta}&P\ar[d]^{\gamma}\\
C\ar[r]^{a}&X_0}$$
is a weak pushout. Thus there exists a morphism $h\colon X_0\to B$ which makes the following
diagram commutative.
$$\xymatrix{
\Omega X_1 \ar[r]^u \ar[d]_{\beta} &P \ar[d]^{\gamma} \ar@/^/[ddr]^w\\
C \ar[r]^a \ar@/_/[drr]_c &X_0 \ar@{-->}[dr]^h\\
&&B}
$$
which implies $\overline{c}=0$. Hence $\overline{\beta}$ is an epimorphism in $\A$.  \qed
\medskip

The following lemma can be found in \cite[Proposition 1.20]{ln}.
\begin{lemma}\label{y1}
Let $\C$ be an extriangulated category and $A\overset{f}{\longrightarrow}B\overset{g}{\longrightarrow}C\overset{\delta}{\dashrightarrow}$
be any $\E$-triangle in $\C$. Assume that $x\colon A\to D$ is any morphism in $\C$. Then there exists
a commutative diagram
$$\xymatrix{
A \ar[r]^f \ar[d]_{x} & B\ar[r]^g \ar[d]^{y} & C\ar@{-->}[r]^{\delta} \ar@{=}[d]&\\
D\ar[r]^{a} & F \ar[r]^{b} &C \ar@{-->}[r]^{x_{\ast}\delta}&}$$
of $\E$-triangles in $\C$, and moreover
$$A\xrightarrow{~\binom{f}{x}~}B\oplus D\xrightarrow{~(y,~-a)~}F\xymatrix{\ar@{-->}[r]^{b^{\ast}\delta}&}$$
becomes an $\E$-triangle in $\C$.
\end{lemma}

\begin{lemma}\label{y2}
Let $\C$ be an extriangulated category with enough projectives  and enough injectives.
Suppose that $\X$ is a cluster-tilting subcategory of $\C$ and $\A$ is the abelian quotient category $\C/\X$.

\emph{(1)} If $f\colon A\to B$ is a morphism in $\C$, then there exists an inflation $\alpha=\binom{f}{a}\colon A\to X_0\oplus B$ in $\B$ such that $\overline {\alpha}=\overline f$.

\emph{(2)} If $f\colon A\to B$ is a morphism in $\C$,
then there exists a deflation $\beta=(f,-b)\colon X_1\oplus A\to B$ in $\B$ such that $\overline \beta=\overline f$.
\end{lemma}

\proof We only show the first one, the second is dual.
Since $\X$ is cluster tilting, there exists an $\E$-triangle
$$A\xrightarrow{~a~}X_0\xrightarrow{~b~}X_1\overset{}{\dashrightarrow},$$
where $X_0,X_1\in\X$. By Lemma \ref{y1}, we get following commutative diagram made of $\E$-triangles.
$$\xymatrix{
A \ar[r]^a\ar[d]_{f} & X_0\ar[r]^b \ar[d]^{y} & X_1\ar@{-->}[r] \ar@{=}[d]&\\
B\ar[r]^{c} & C \ar[r]^{d} &X_1 \ar@{-->}[r]&}$$
Moreover
$A\xrightarrow{~\alpha=\binom{f}{a}~}X_0\oplus B\xrightarrow{~(y,~-c)~}C\overset{}{\dashrightarrow}$
is an $\E$-triangle in $\C$.

This shows that $\alpha$ is an inflation and $\overline {\alpha}=\overline f$.  \qed

\begin{lemma}\label{x1}
Let $\C$ be an extriangulated category with enough projectives and enough injectives.
Suppose that $\X$ is a cluster-tilting subcategory of $\C$. Then
$\Omega\X$ and $\Sigma\X$ are closed under direct summands.
\end{lemma}

\proof See the proof of Lemma 5.9 in \cite{ln}.   \qed

\begin{remark}\label{rem1}
Let $\C$ be an extriangulated category with enough projectives and enough injectives.
If $\X$ is a cluster-tilting subcategory of $\C$, then
$\Omega\X/\X=\Omega\X/\cal P$ and $\Sigma\X/\X=\Sigma\X/\cal I$. For convenience, we denote
$\Omega\overline{\X}:=\Omega\X/\X$ and $\Sigma\underline{\X}:=\Sigma\X/\X$.
\end{remark}

\proof We only prove $\Omega\X/\X=\Omega\X/\cal P$. By duality, we have $\Sigma\X/\X=\Sigma\X/\cal I$.

We first prove that a morphism $f\colon \Omega X\to C$ factors through $\cal P$ with $X\in\X$  if and only if it factors through $\X$. Since $\cal P\subseteq\X$, we have $\cal P\subseteq\X$. We only need to prove
$f$ factors through $\X$ implies it factors through $\cal P$. Assume that $f$ factors through $\X$.
By the definition of $\Omega\X$, we have the following $\E$-triangle:
$$\Omega X\overset{a}{\longrightarrow}P\overset{b}{\longrightarrow}X\overset{}{\dashrightarrow},$$
where $P\in\cal P$. Since $\X$ is cluster tilting, there exists an
$\E$-triangle:
$$X_0\overset{c}{\longrightarrow}X_1\overset{d}{\longrightarrow}C\overset{}{\dashrightarrow},$$
where $X_0,X_1\in\X$.
Since $\E(\X,\X)=0$, we have that $d$ is a right $\X$-approximation of $C$. Then there exists a morphism
$g\colon \Omega X\to X_1$ such that $f=dg$. Since $a$ is a left $\X$-approximation of $\Omega X$, there exists a morphism $h\colon P\to X_1$ such that $g=ha$. It follows that $f=(dh)a$. This shows that $f$ factors through $\cal P$.

Thus by definition we have $\Omega\X/\X=\Omega\X/\cal P$.  \qed

\begin{lemma}\label{y3}
Let $\C$ be an extriangulated category with enough projectives  and enough injectives.
Suppose that $\X$ is a cluster-tilting subcategory of $\C$ and $\A$ is the abelian quotient category $\C/\X$.
Then an object $M$ of $\A$ is a projective object if and only if $M\in \Omega\overline{\X}$. Dually an object $N$ of $\A$ is an injective object if and only if $N\in\Sigma\underline{\X}$.
\end{lemma}

\proof We prove the first statement only, the second one is obtained dually.

Let $\overline{g}\colon B\to C$ be an epimorphism in $\A$ and $\overline{\beta}\colon \Omega X\to C$ be any morphism
in $\C$ where $X\in\X$. By Lemma \ref{y2}, we can assume that it admits an $\E$-triangle
$$A\xrightarrow{~f~}B\xrightarrow{~g~}C\overset{}{\dashrightarrow}$$
Since $\X$ is cluster tilting, there exists an $\E$-triangle
$$B\xrightarrow{~a~}X_0\xrightarrow{~b~}X_1\overset{}{\dashrightarrow},$$
where $X_0,X_1\in\X$. By (ET4), we get following commutative diagram made of $\E$-triangles.
\begin{equation}\label{t1}
\begin{array}{l}
\xymatrix{A\ar[r]^{f}\ar@{=}[d]&B\ar[r]^{g}\ar[d]^{a}&C\ar[d]^{u}\\
A\ar[r]^{c}&X_0\ar[d]^{b}\ar[r]^{d}&D\ar[d]^{v}\\
&X_1\ar@{=}[r]&X_1}
\end{array}
\end{equation}
It follows that $ug=da$ and then $\overline{u}\circ\overline{g}=0$. Since $\overline{g}$
is an epimorphism, we have $\overline{u}=0$.
By definition $\Omega X$ admits an $\E$-triangle
$\Omega X\xrightarrow{~p~}P\xrightarrow{~q~}X\overset{}{\dashrightarrow}$
where $P\in\cal P$.
Since $\overline{u}\circ \overline{\beta}=0$, then $u\beta$ factors through $\X$.
As $\E(\X,\X)=0$,  we obtain that $p$ is a left $\X$-approximation of $\Omega X$.
Thus there exists a morphism $r\colon P\to D$ such that
$rp=u\beta$.  Since $P$ is a project object, there exists a morphism
$w\colon P\to X_0$ such that $r=dw$. It follows that $d(wp)=u\beta$.
By the dual of \cite[Lemma 3.13]{np}, the upper-right square in the diagram (\ref{t1})
$$\xymatrix{B\ar[r]^{g}\ar[d]^{a}&C\ar[d]^{u}\\
X_0\ar[r]^{d}&D}$$
is a weak pullback. Thus there exists a morphism $h\colon \Omega X\to B$ which makes the following
diagram commutative.
$$\xymatrix{
\Omega X \ar@/^/[drr]^{\beta} \ar@{-->}[dr]^h \ar@/_/[ddr]_{wp}\\
&B \ar[r]^g \ar[d]^a &C \ar[d]^u\\
&X_0 \ar[r]^d &D}
$$
Hence $\overline{\beta}=\overline{g}\circ\overline{h}$. This shows that $\Omega X$ is a projective object in $\A$.
\smallskip

Conversely, assume that $M$ is a projective object in $\A$, by Lemma \ref{x2}, there exists
an epimorphism $\overline{\beta}\colon \Omega X\to M$ for some $X\in\X$ in $\A$.
Thus $M$ is a direct summand of $\Omega X$ in $\A$. Hence by Lemma \ref{x1}, we have $M$ lies in $\Omega\overline{\X}$.  \qed

\medskip
Recall that an abelian category with enough projectives and injectives is called \emph{Gorenstein}
if all projective objects of this category have finite injective dimension, and all injective
objects have finite projective dimension. The maximum of the injective dimensions of projectives and the projective dimensions of injectives is called \emph{Gorenstein dimension} of the category.

\begin{theorem}\label{thm}
Let $\C$ be an extriangulated category with enough projective objects  and enough injective objects.
Suppose that $\X$ is a cluster-tilting subcategory of $\C$ and $\A$ is the abelian quotient category $\C/\X$. Then:
\begin{itemize}
\item[\emph{(1)}] The category $\A$ has enough projective objects and enough injective objects.
\item[\emph{(2)}] If $\Sigma(\Omega\X)\subseteq\X$ and $\Omega(\Sigma\X)\subseteq\X$, then the category $\A$ is Gorenstein of Gorenstein dimension at most one.
\end{itemize}
\end{theorem}

\proof (1) This follows from Lemma \ref{x2} and Lemma \ref{y3}.

(2) Let $\Sigma X$ be any injective object in $\A$. Since $\X$ is cluster tilting, there exists an $\E$-triangle
$$\Sigma X\xrightarrow{~a~}X_0\xrightarrow{~b~}X_1\overset{}{\dashrightarrow},$$
where $X_0,X_1\in\X$. By the definition of $\Omega\X$, we have the following $\E$-triangle:
$$\Omega X_0\overset{u}{\longrightarrow}P\overset{v}{\longrightarrow}X_0\overset{}{\dashrightarrow},$$
where $P_0\in\cal P$.  By (ET4)$\op$, we have the following commutative diagram made of $\E$-triangles
\begin{equation}\label{t3}
\begin{array}{l}
\xymatrix{\Omega X_0\ar[d]^p\ar@{=}[r]&\Omega X_0\ar[d]^u&&\\
\Omega X_1\ar[r]^c\ar[d]^{q}&P_0\ar[r]^d\ar[d]^{v}&X_1\ar@{=}[d]\ar@{-->}[r]&\\
\Sigma X\ar[r]^{a}\ar@{-->}[d]&X_0\ar[r]^{b}\ar@{-->}[d]&X_1\ar@{-->}[r]&\\
&&&}
\end{array}
\end{equation}
By the definition of $\Omega\X$, we have the following $\E$-triangle:
$$\Omega(\Sigma X)\overset{x}{\longrightarrow}P_1\overset{y}{\longrightarrow}\Sigma X\overset{}{\dashrightarrow},$$
where $P_1\in\cal P$.
By the dual of \cite[Proposition 3.17]{np}, we obtain the following commutative
diagram made of $\E$-triangles.
\begin{equation}\label{t4}
\begin{array}{l}
\xymatrix@C=1.2cm@R=1.2cm{&\Omega X_0\ar@{=}[r]\ar[d]^{\binom{0}{1}}&\Omega X_0\ar[d]^{p}\\
\Omega(\Sigma X)\ar[r]^{\binom{x}{h}\quad\;}\ar@{=}[d]&P_1\oplus \Omega X_0\ar[r]^{\quad (-p',\ p)}\ar[d]^{(1,\ 0)}&\Omega X_1\ar@{-->}[r]\ar[d]^{q}&\\
\Omega(\Sigma X)\ar[r]^{x}&P_1\ar[r]^{y}\ar@{-->}[d]&\Sigma X\ar@{-->}[r]\ar@{-->}[d]&\\
&&}
\end{array}
\end{equation}
We claim that
$$\Omega(\Sigma X)\xrightarrow{~\overline{h}~}\Omega X_0\xrightarrow{~\overline{p}~}\Omega X_1\xrightarrow{~\overline{q}~}\Sigma X\xrightarrow{~~} 0$$
is an exact sequence in $\A$. In fact, in the diagram (\ref{t4})
 we obtain that $qp=0$ and $$(-p',p)\binom{x}{h}=0$$ which implies
$\overline{q}\circ \overline{p}=0$ and $\overline{p}\circ\overline{h}=0$.
This shows that ${\rm Im}(\overline{p})\subseteq {\rm Ker}(\overline{q})$ and ${\rm Im}(\overline{h})\subseteq {\rm Ker}(\overline{p})$.

Now we show that ${\rm Ker}(\overline{q})\subseteq{\rm Im}(\overline{p})$.

Let $\overline{\alpha}\colon M\to \Sigma X$ be any morphism in $\A$ such that $\overline{q}\circ\overline{\alpha}=0$. Then $q\alpha$ factors through $\X$. By Remark \ref{rem1}, we know that
$q\alpha$ factors through $\cal P$. That is to say, there exist morphisms $s\colon M\to P_2 $
and $t\colon P_2\to \Sigma X$ such that $q\alpha=ts$ where $P_2\in \cal P$.
Since $P_2$ is a project object, there exists a morphism $\beta\colon P_2\to \Omega X_1$ such that
$q\beta=t$ and then $$q(\alpha-\beta s)=q\alpha-q\beta s=q\alpha-ts=0.$$
Thus there exists a morphism $\gamma\colon M\to \Omega X_0$ such that $\alpha-\beta s=p\gamma$
and then $\alpha=\beta s+p\gamma$.
It follows that $\overline{\alpha}=\overline{p}\circ\overline{\gamma}$ which implies ${\rm Ker}(\overline{q})\subseteq{\rm Im}(\overline{p})$.

Now we show that ${\rm Ker}(\overline{p})\subseteq {\rm Im}(\overline{h})$.

Let $\overline{l}\colon N\to \Omega X_0$ be any morphism in $\A$ such that $\overline{p}\circ\overline{l}=0$.
Then $pl$ factors through $\X$. By Remark \ref{rem1}, we know that
$pl$ factors through $\cal P$. That is to say, there exist morphisms $f\colon N\to P_3 $
and $g\colon P_3\to \Omega X_1$ such that $pl=gf$ where $P_3\in \cal P$.
Since $P_3$ is a project object, there exists a morphism $\binom{m}{n}\colon P_3\to P_1\oplus\Omega X_0$ such that
$g=(-p',p)\binom{m}{n}=-p'm+pn$ and then
$$(-p',p)\binom{mf}{nf-l}=(-p'm+pn)f-pl=0.$$
Thus there exists a morphism $w\colon N\to \Omega(\Sigma X)$ such that $\binom{x}{h}w=\binom{mf}{nf-l}$
and then $l=nf-hw$.
It follows that $\overline{l}=\overline{h}\circ(-\overline{w})$ which implies ${\rm Ker}(\overline{p})\subseteq{\rm Im}(\overline{h})$.

Now we show that $\overline{q}$ is an epimorphism in $\A$.

Let $\overline{i}\colon \Sigma X\to L$ is any morphism in $\A$ such that
$\overline{i}\circ\overline{q}=0$. Then $iq$ factors through $\X$, namely, there exist morphisms
$i\colon \Omega X_1\to X_2$ and $k\colon X_2\to L$ such that $iq=kj$.
Since $\E(\X,\X)=0$, we have that $c$
is a left $\X$-approximation of $\Omega X_1$. Thus there exists a morphism
$k'\colon P_0\to X_2$ such that $k'c=j$. It follows that $iq=(kk')c$.
By \cite[Lemma 3.13]{np}, the lower-left square in the diagram (\ref{t3})
$$\xymatrix{\Omega X_1\ar[r]^{c}\ar[d]^{q}&P_0\ar[d]^{v}\\
\Sigma X\ar[r]^{a}&X_0}$$
is a weak pushout. Thus there exists a morphism $z\colon X_0\to L$ which makes the following
diagram commutative.
$$\xymatrix{
\Omega X_1 \ar[r]^c \ar[d]_{q} &P_0 \ar[d]^{v} \ar@/^/[ddr]^{kk'}\\
\Sigma X \ar[r]^a \ar@/_/[drr]_i &X_0 \ar@{-->}[dr]^z\\
&&B}
$$
which implies $\overline{z}=0$. Hence $\overline{q}$ is an epimorphism in $\A$.
\medskip

This shows that $\Omega(\Sigma X)\xrightarrow{~\overline{h}~}\Omega X_0\xrightarrow{~\overline{p}~}\Omega X_1\xrightarrow{~\overline{q}~}\Sigma X\xrightarrow{~~} 0$
is an exact sequence in $\A$.

In the diagram (\ref{t3}), we obtain that $aq=vc$ and $u=cp$. In the  diagram (\ref{t4}), we obtain that $y=-qp'$ and $p'x=ph$.
Thus we have that $ay=-aqp'=v(-cp')$ and $$(-cp')x=-cph=-uh.$$ Hence we have the following commutative diagram of $\E$-triangles.
$$\xymatrix{
\Omega(\Sigma X)\ar[r]^{\quad x}\ar[d]_{-h} & P_1\ar[r]^y \ar[d]^{-cp'} &\Sigma X\ar@{-->}[r] \ar[d]^{a}&\\
\Omega X_0\ar[r]^{u} &P_0 \ar[r]^{v} &X_0 \ar@{-->}[r]&}$$
By the definition of $\Omega$, we have $\Omega a=-h$ and then $\overline{h}=-\Omega\overline{a}$.
Since $\Omega a\colon \Omega(\Sigma X)\to \Omega X_0$ and $\Omega(\Sigma\X)\subseteq\X$, we have
$\Omega \overline{a}=0$ in $\A$. Namely $\overline{h}=0$ in $\A$.
So we obtain that
$$0\xrightarrow{~~}\Omega X_0\xrightarrow{~\overline{p}~}\Omega X_1\xrightarrow{~\overline{q}~}\Sigma X\xrightarrow{~~} 0$$
is an exact sequence in $\A$.

This shows that any injective object $\Sigma X$ in $\A$ has projective dimension at most one.

Dually, we can show that any projective object in $\A$ has injective dimension at most one.

Therefore $\A$ is Gorenstein of Gorenstein dimension at most one. \qed
\medskip

We conclude this section with an example illustrating our result.

\begin{example}
Let $\Lambda$ be a finite dimensional algebra given by the quiver
$$\xymatrix@C1cm{1 &2\ar[l]_{\;\;\;\alpha}&3\ar[l]_{\;\;\;\beta} }$$
with relation $\beta\alpha=0$. Let $S(i)$ be the simple module concentrated at the vertex $i$, and $P(i)$ be the indecomposable projective right $\Lambda$-module.  The
AR-quiver of ${\rm mod}\Lambda$ is given by
 $$\xymatrix@C=0.4cm@R0.4cm{&P(2)\ar[dr]\ar@{.}[rr]&&P(3)\ar[dr]\\
 P(1)\ar[ur]\ar@{.}[rr]&&S(2)\ar[ur]\ar@{.}[rr]&&S(3)}$$
It is straightforward to verify that the subcategory
$$\X=\add\big(P(1)\oplus P(2)\oplus P(3)\oplus S(3)\big)$$
is a cluster tilting subcategory (is also called maximal $1$-orthogonal subcategory) of ${\rm mod}\Lambda$.
Since $\Sigma(\Omega\X)=\add\big(S(3)\big)\subseteq\X$ and $\Omega(\Sigma\X)=\add\big(P(1)\big)\subseteq\X$.
By Theorem \ref{thm}, we have that ${\rm mod}\Lambda/\X$ is  Gorenstein of Gorenstein dimension at most one.
\end{example}

Yu Liu\\
School of Mathematics, Southwest Jiaotong University, 610031, Chengdu,
Sichuan, P. R. China\\
E-mail: \textsf{liuyu86@swjtu.edu.cn}\\[0.3cm]
Panyue Zhou\\
College of Mathematics, Hunan Institute of Science and Technology, 414006, Yueyang, Hunan, P. R. China.\\
E-mail: \textsf{panyuezhou@163.com}


\begin{thebibliography}{99}
\addtolength{\itemsep}{-0.5em}
\bibitem[DL]{dl} L. Demonet, Y. Liu. Quotients of exact categories by cluster tilting subcategories as module categories.  J. Pure Appl. Algebra 217(12),  2282-2297, 2013.


\bibitem[KZ]{kz}
S. Koenig, B. Zhu.
From triangulated categories to abelian categories: cluster tilting in a general framework.
 Math. Z. 258, 143-160, 2008.



\bibitem[Liu]{liu} Y. Liu. Abelian quotients associated with fully rigid subcategories. arXiv: 1902.07421, 2019.


\bibitem[LN]{ln} Y. Liu, H. Nakaoka. Hearts of twin cotorsion pairs on extriangulated categories. J. Algebra 528: 96-149, 2019.


\bibitem[NP]{np} H. Nakaoka, Y. Palu. Mutation via Hovey twin cotorsion pairs and model structures in extriangulated categories.  arXiv: 1605.05607, 2016.

\bibitem[ZZ]{zz} P. Zhou, B. Zhu. Cluster-tilting subcategories in extriangulated categories. Theory Appl. Categ.
 34(8): 221-242, 2019.
 
\bibitem[ZZ1]{zz1} P. Zhou, B. Zhu. Triangulated quotient categories revisited. J. Algebra 502: 196-232, 2018.




\end{thebibliography}
\end{document}